\documentclass[preprint,11pt]{elsarticle}
\usepackage{amssymb}
\usepackage{mathrsfs}
\usepackage{amsthm}
\usepackage{amsmath}
\usepackage{epstopdf}
\usepackage{epsfig}
\usepackage{color}
\usepackage{graphicx}
\usepackage{subfigure}
\usepackage{enumerate}
\usepackage{longtable,bookt
abs}
\usepackage{lineno,hyperref}
\usepackage{threeparttable}
\usepackage{cleveref}
\usepackage{verbatim}
\usepackage{comment}

\newtheorem{theorem}{Theorem}[section]
\newtheorem{lemma}{Lemma}[section]
\newtheorem{corollary}{Corollary}[section]
\newtheorem{proposition}{Proposition}[section]
\newtheorem{definition}{Definition}[section]
\newtheorem{remark}{Remark}[section]

\begin{document}
\begin{frontmatter}


\title{Permanence for continuous-time competitive Kolmogorov systems via the carrying simplex}

\author[add1,add2]{Lei Niu}\ead{lei.niu@dhu.edu.cn}
\author[add1]{Yuheng Song}
\address[add1]{Department of Mathematics, Donghua University, Shanghai
201620, China}
\address[add2]{Institute for Nonlinear Science, Donghua University, Shanghai 201620, China}
\begin{abstract}
In this paper we study the permanence and impermanence for continuous-time competitive Kolmogorov systems via the carrying simplex. We first give an extension to attractors of V. Hutson's results on the existence of repellors in continuous-time dynamical systems  that have found wide use in the study of permanence via average Liapunov functions. We then give a general criterion for the stability of the boundary of carrying simplex for competitive Kolmogorov systems of differential equations, which determines the permanence and impermanence of such systems. Based on the criterion, we present a complete classification of the permanence and impermanence in terms of inequalities on parameters for all three-dimensional competitive systems which have linearly determined nullclines. The results are applied to many classical models in population dynamics including the Lotka-Volterrra system, Gompertz system, Leslie-Gower system and Ricker system.

\end{abstract}


\end{frontmatter}



\section{Introduction}
Permanence is an important property of dynamical systems, and of the ecological systems in ecology and biology they are modeling. The theory addresses under which conditions a collection of interacting species can coexist over long periods of time, which has led to the development of many interesting mathematical techniques and results. The mathematical investigation of permanence (also known as uniform persistence) began with the early work of Freedman and Waltman \cite{Freedman1977}, Gard \cite{Gard1980}, Gard and Hallam \cite{Gard1979}, Schuster Sigmund and Wolff \cite{Schuster1979}, among others, in the late seventies. Since then, the theory has developed rapidly using the available tools from dynamical system theory, which simply requires that the set of all extinction states, be a repellor for the systems (see \cite{Butler1986,Freedman1989,Garay1989,Gyllenberg2020b,Hofbauer1989,Hutson1984,Smith2011,Zhao2017}). One of the most famous results is Hutson’s results in \cite{Hutson1984} on the existence of repellors via the average Liapunov functions, which have found wide
use in the study of permanence for the Kolmogorov systems (see, for example, \cite{Garay2003,hofbauer1987coexistence,Hutson1988,Kon2004,Lu2017, Lu1999}). The average Liapunov function method is first used by Hofbauer in \cite{Hofbauer1981} and further generalized by Hutson and Moran \cite{Hutson1982} for discrete-time dynamical systems and by Hutson \cite{Hutson1984} for continuous-time dynamical systems.

In this paper, we study the permanence and impermanence for the continuous-time Kolmogorov systems of differential equations
\begin{equation}\label{kolmogorov-equas}
\dot{x}_i=x_if_i(x),\quad x_i\geq 0,\quad i=1,\ldots,n
\end{equation}
when it admits a carrying simplex, that is a globally attracting invariant hypersurface of codimension one. Due to Hirsch's work in \cite{hirsch1988}, the continuous-time competitive Kolmogorov systems of type \eqref{kolmogorov-equas} can possess a carrying simplex. A weaker version of Hirsch's result was provided recently by Hou in \cite{hou2020}. The theory of the carrying simplex has also been extended to competitive Kolmogorov mappings  (see \cite{diekmann2008carrying, hirsch1988, hou2021, LG, jiang2015,ruiz2011exclusion,wang2002uniqueness}), since the early work of de Mottoni and Schiaffino \cite{Mottoni1981}, Hale and Somolinos \cite{Hale1983} and Smith \cite{smith1986}. The importance of the existence of a carrying simplex stems from the fact that it captures all the relevant long-term dynamics. Our project here is to provide the minimal conditions to ensure the permanence and impermanence for such systems via the carrying simplex. The main mathematical tools involve average Liapunov functions and the theory of the carrying simplex. Successful cases of combining these two approaches include the stability criteria of heteroclinic cycles provided in \cite{jiang2015} and permanence criteria in \cite{Gyllenberg2020b} for competitive Kolmogorov mappings. 

Inspired by the methods in \cite{jiang2015}, we first provide an extension (Theorem \ref{attractor_T}) to attractors of Hutson's results in \cite{Hutson1984} on the existence of repellors in continuous-time dynamical systems. This result could be perceived not only as a technique for studying the impermanence of the dynamical systems but has its own interest. We then extend the topological results on repellors and attractors to dissipative Kolmogorov systems (Theorem \ref{theorem:permanence}), and we give a general criterion (Theorem \ref{Sigma-permanence}) for the stability of the boundary of the carrying simplex for continuous-time competitive Kolmogorov systems \eqref{kolmogorov-equas}, which can determine the permanence and impermanence of the systems. In particular, for three-dimensional competitive systems, permanence and impermanence can be determined by boundary equilibria (Corollary \ref{coro:3-permanence}). As an interesting case, when the boundary of the carrying simplex is a heteroclinic cycle, the system is permanent if the heteroclinic cycle is a repellor, while the system is impermanent if it is an attractor (Corollary \ref{coro:h_cycle}).

By using our criteria, we give a detailed analysis of the three-dimensional competitive Kolmogorov systems
\begin{equation}\label{equ:cssn}
    \frac{dx_i}{dt}=x_if_i(c_i,\sum_{j=1}^3b_{ij}x_j)=x_if_i(c_i,(Bx^\tau)_i), \quad i=1,2,3
\end{equation}
 with linearly determined nullclines which are often used to describe population dynamics  (see, for example, \cite{franke1991,hou2019,JN2017,levine2002,rees1997}), 
where $c_i,b_{ij}>0$, $B=(b_{ij})_{3\times 3}$, and $f_i:\dot{\mathbb{R}}_+\times \dot{\mathbb{R}}_+\to \mathbb{R}$ are $C^1$ functions such that
\begin{equation}\label{property:f}
    \left\{\begin{array}{l}
                 \displaystyle\lim_{y\to 0}yf_i(r,y)=0, \\
                 \noalign{\medskip}

                 \frac{\partial f_i(r,y)}{\partial y}<0 ~\textrm{and} ~f_i(r,r)=0,\quad \forall (r,y)\in \dot{\mathbb{R}}_+\times \dot{\mathbb{R}}_+.
 \end{array}\right.
\end{equation}
The classical Lotka-Volterrra system (\cite{Zeeman1994,Z993}), Gompertz system (\cite{jnz,Lu2017}), Leslie-Gower system (\cite{hou2019, JN2017}) and Ricker system (\cite{hou2019, JN2017}) are specific examples of such systems of type \eqref{equ:cssn}. System \eqref{equ:cssn} possesses a carrying simplex, and Jiang and Niu \cite{JN2017} showed that the dynamics of \eqref{equ:cssn} can always be classified into 33 classes by nullcline equivalence, where the classification method was first developed by Zeeman for three-dimensional competitive Lotka-Volterrra system in \cite{Z993}. Here, we illustrate how to utilize our criteria to study the  permanence and impermanence for competitive systems \eqref{equ:cssn}. In particular, we can give a complete classification of the permanence and impermanence in terms of inequalities on parameters for \eqref{equ:cssn} based on Jiang and Niu's classification in \cite{JN2017}. Specifically, we prove that system \eqref{equ:cssn} is permanent if it is in Jiang and Niu's classes 29, 31, 33 and class 27 with a repelling heteroclinic cycle, while it is impermanent if it is in their classes 1--26, 28, 30, 32 and class 27 with an attracting heteroclinic cycle (Therorem \ref{permanence-3D}). 

The paper is organized as follows. Section \ref{sec:2} includes some notation. In Section \ref{sec:3}, we provide an extension to attractors of Hutson's results in \cite{Hutson1984} on the existence of repellors in continuous dynamical systems. In Section \ref{sec:4}, we extend the topological results on repellors and attractors to dissipative Kolmogorov systems and competitive systems and present the criteria on permanence and impermanence. We then give a detailed analysis for three-dimensional competitive systems. In Section \ref{sec:5}, we apply our results to some classical population
models of three competing species, and give the completed classification of permanence and impermanence  via the equivalence classes provided by Jiang and Niu in \cite{JN2017}. 

\section{Notation}\label{sec:2}
Let $X$ be a metric space with metric $d_X(\cdot,\cdot)$ and let $\mathbb{Z}_+=\{0,1,2,\ldots\}$, $\mathbb{R}_+=[0,+\infty)$. Suppose that $(X, \mathbb{R}_+, \varPsi)$ is a continuous semi-dynamical system with the state space $X$ and the semi-flow $\varPsi$, that is, $\varPsi: \mathbb{R}_+\times X \to X$, $(t,x)\mapsto \varPsi_t(x)$ is a continuous map such that 
$\varPsi_0(x)=x$ and $\varPsi_t(\varPsi_s(x))=\varPsi_{t+s}(x)$ for all $t,s\geq 0$, $x\in X$. The (positive) orbit of $\varPsi$ through $x\in X$ is denoted by $\gamma^+(x,\varPsi):=\{\varPsi_t(x): t\geq 0\}$, and by $\gamma^+(x)$ when $\varPsi$ is understood. Denote the tail from the moment $\tau>0$ of $\gamma^+(x)$ by $\gamma^+_\tau(x):=\{\varPsi_t(x): t\geq \tau\}$. When $\gamma^+(x) = \{x\}$ then $x$ is an equilibrium (or fixed point). The set of equilibria is denoted by $\mathcal{E}(\varPsi)$. The omega limit set $\omega(x)$ of $x$ is the set of limit points of the positive orbit $\gamma^+(x)$.

For a subset $D\subseteq X$, $\overline{D}$ and $D^c$ denote the closure of $D$ in $X$ and the complement of $D$ respectively. Given any set $D$, let $\gamma^+(D)=\bigcup_{x\in D}\gamma^+(x)$ be the positive orbit through $D$ and define the omega limit set $\omega(D)$ of $D$ as $\omega(D)=\bigcap_{s \geq 0} \overline{\bigcup_{t \geq s} \varPsi_s(D)}$. Let $\Lambda(D)=\bigcup_{x\in D}\omega(x)$. 

 A set $D\subseteq X$ is positively invariant (under $\varPsi$) if $\varPsi_t(D)\subseteq D$ for all $t\geq 0$; negatively invariant if $\varPsi_t(D)\supseteq D$ for all $t\geq 0$; and invariant if $\varPsi_t(D)=D$ for all $t\geq 0$. The set $B$ is said to be absorbing for set $U$ if it is positively invariant and $\gamma^+(x)\cap B\neq \emptyset$ for every $x\in U$.
 
We say that $J\subseteq X$ attracts a set $D$ under $\varPsi$ if for any $\varepsilon > 0$, there exists a $t_0\geq 0$ such that $\varPsi_t(D)$ belongs to the $\varepsilon$-neighborhood $O_\varepsilon(J)$ of $J$ for $t > t_0$, where $O_\varepsilon(J)=\{y\in X:d_X(y,J)<\varepsilon\}$; when $D =\{x\}$ we say $J$ attracts $x$. A  subset $J$ is said to be stable under $\varPsi$ if for any neighborhood $\mathcal{V}$ of $J$, there is a neighborhood $\mathcal{U}$ of $J$ such that $\varPsi_t(\mathcal{U}) \subseteq \mathcal{V}$ for all $t \geq 0$. A nonempty compact invariant set $J\subseteq X$ is said to be an attractor of $\varPsi$ if it attracts a neighborhood of itself; $J$ is called a global attractor of $\varPsi$ if it attracts each bounded set $D \subseteq X$ (see \cite{Akin1993,Conley1978,Hale1991}). If $\varPsi$ has a global attractor then it is said to be dissipative. Suppose $S\subseteq X$ is a nonempty, compact and positively invariant subset. We say that $S$ is a repellor (under $\varPsi$) if there exists  a neighborhood $\mathcal{U}$ of $S$ such that for all $x \notin S$ there exists a $t_0=t_0(x)>0$ satisfying $\varPsi_t(x)\notin \mathcal{U}$ for all $t \geq t_0$ (see \cite{Hirsch2001,Hofbauer1989}).

Throughout this paper, we reserve the symbol $n$ for the dimension of the Euclidean space $\mathbb{R}^n$ and the symbol $N$ for the set $\{1, \ldots, n\}$. We will denote by $\{e_{\{1\}}, \ldots, e_{\{n\}}\}$ the usual basis for $\mathbb{R}^n$, and by $d(\cdot,\cdot)$ the usual Euclidean distance. We use $\mathbb{R}^n_+$ to denote the nonnegative cone $\{x\in \mathbb{R}^n: x_i\geq 0,\forall i\in N\}$. The interior of $\mathbb{R}^n_+$ is the open cone $\dot{\mathbb{R}}^n_+:= \{x\in \mathbb{R}^n_+: x_i>0,\forall i\in N \}$ and the boundary of $\mathbb{R}^n_+$ is $\partial \mathbb{R}^n_+:=\mathbb{R}^n_+\setminus \dot{\mathbb{R}}^n_+$. We denote by $\mathbb{H}^+_{\{i\}}$ the $i$-th positive coordinate axis and by $\pi_i=\{x\in \mathbb{R}^n_+:x_i=0\}$ the $i$-th coordinate plane. For $x\in \mathbb{R}^n_+$, let $\kappa(x)=\{i:x_i>0\}$ denote its support and $x^\tau$ denote its transpose. The symbol $0$  stands for both the origin of $\mathbb{R}^n$ and the real number $0$.

\section{Topological results on repellor and attractor}\label{sec:3}
Throughout this section, $M$ will be a compact metric space with metric $d_M(\cdot,\cdot)$ and $S$ is a compact subset of
 $M$ with empty interior. Suppose that $\varPsi:\mathbb{R}_+\times M\to M$ is a continuous semi-flow. 

\begin{lemma}[Criterion for repellor \cite{Hutson1984}]\label{repellor_T}
Let $M$ be a compact metric space and consider the continuous dynamical system $(M,\mathbb{R}_+,\varPsi)$. Assume that $S$ is a compact subset of $M$ with empty interior such that $S$ and $M\setminus S$ are positively invariant under $\varPsi$. Suppose that there is a continuous function $V: M\to \mathbb{R}_+$ satisfying that
\smallskip

  {\rm (i)}~$V(x)=0\Leftrightarrow x\in S$,
\smallskip

  {\rm (ii)}~$\displaystyle\psi_M(x):=\sup_{t\geq 0}\vartheta_M(t,x)>1$ for all $x\in S$, where
 \begin{equation}
 \vartheta_M(t,x)=\liminf_{\mbox{\tiny$\begin{array}{c}
y\rightarrow x\\
y\in M\setminus S
\end{array}$}}\!\frac{V(\varPsi_t(y))}{V(y)},~~t\geq 0.
 \end{equation} 
Then $S$ is a repellor.
\end{lemma}
\begin{remark}\label{con:repel}
The function $V$ in Lemma \ref{repellor_T} is called an average Liapunov function. Moreover, by \cite[Corollary 2.3]{Hutson1984}, the condition {\rm (ii)} in Lemma \ref{repellor_T} is implied by the following condition

\smallskip

{\rm (ii')}~$\psi_M(x)>1$ for all $x\in \Lambda(S)$, and $\psi_M(x)>0$ for all $x\in S$. 
\smallskip

\noindent By following the proof of \cite[Proposition 2.4]{Gyllenberg2020b} and \cite[Lemma 2.1]{Hutson1984}, it is easy to show that if $S$ is repellor, then there is a compact positively invariant set $K\subseteq M\setminus S$ such that for every $x\in M\setminus S$, there exists a $m=m(x)>0$ satisfying $\gamma^+_m(x)\subseteq  K$.
\end{remark}

\begin{theorem}[Criterion for attractor]\label{attractor_T}
Consider the dynamical system $(M,\mathbb{R}_+,\varPsi)$ as in Lemma \ref{repellor_T}.  Suppose that there is a continuous function $V: M\to \mathbb{R}_+$ and some $t_0>0$, $C>1$ such that
\smallskip

{\rm (i)}  $V(x)=0\Leftrightarrow x\in S$, and $\frac{V(\varPsi_t(x))}{V(x)}\leq C$ for all $x\in M\setminus S$ and $0\leq t\leq t_0$,
\smallskip

{\rm (ii)} $\displaystyle\varphi_M(x):=\inf_{t\geq 0}\zeta_M(t,x)<1$ for all $x\in \Lambda(S)$, where
\begin{equation}
\zeta_M(t,x)=\limsup_{\mbox{\tiny$\begin{array}{c}
y\rightarrow x\\
y\in M\setminus S
\end{array}$}}\!\frac{V(\varPsi_t(y))}{V(y)},~~t\geq 0.
\end{equation}
Then $J=\bigcap_{t\geq 0}\varPsi_t(S)$ is an attractor which attracts a neighborhood $\mathcal{U}$ of $S$, and in particular, $S$ is an attractor if it is invariant.
\end{theorem}
\begin{proof}
The proof is similar to that of Theorem 1 in \cite{jiang2015}. For $t \geq 0$, define
$$
   \zeta_M(t,x)= \left \{
    \begin{array}{l}
        V(\varPsi_t(x))/V(x),\quad  x\in M\setminus S, \\[6pt]
        \displaystyle\limsup_{\mbox{\tiny$\begin{array}{c}
y\rightarrow x\\
y\in M\setminus S
\end{array}$}} V(\varPsi_t(y))/V(y), \quad x\in S.
    \end{array}
    \right.
$$
Given $t\geq 0$, there exists an integer $l_t\in \mathbb{Z}_+$ and a real number $0 \leq s<t_0$ such that $t=l_t t_0+s$. By (i) and the positive invariance of $M\setminus S$, for all $t\geq 0$ and $x\in M\setminus S$, one has
$$
    \begin{array}{rl}        \zeta_M(t,x)=&\displaystyle\frac{V(\varPsi_t(x))}{V(x)}\\[6pt]
        \noalign{\medskip}
=&\displaystyle\frac{V(\varPsi_s(\varPsi_{l_t t_0}(x)))}{V(\varPsi_{l_t t_0}(x))}\frac{V(\varPsi_{l_t t_0}(x))}{V(\varPsi_{(l_t-1) t_0}(x))}\frac{V(\varPsi_{(l_t-1) t_0}(x))}{V(\varPsi_{(l_t-2) t_0}(x))}\cdots\frac{V(\varPsi_{t_0}(x))}{V(x)}\\[6pt]
\noalign{\medskip}
\leq & C^{l_t+1}.
    \end{array}
$$
Therefore, $\zeta_M(t,x)\leq C^{l_t+1}<+\infty$ for all $x\in M$, that is, $\zeta_M(t,\cdot):M\to \mathbb{R}_+$ is well defined. Note that for each $t\geq 0$, $\zeta_M(t,\cdot)$ is an upper semi-continuous function, which ensures that 
$U(h,t):=\{x:\zeta_M(t,x)<1-h\}$ 
is an open set for all $0<h<1$ and $t\geq 0$. Let $\varphi_M(x):=\displaystyle\inf_{t\geq 0} \zeta_M(t,x)$. Then by a direct modification of the proof of \cite[Theorem 1]{jiang2015}, one can obtain that $\varphi_M(x)<1$ for all $x\in S$ if (ii) holds, and hence, there exist an $0<h<1$, and $t_1,\ldots,t_m>0$, such that 
 $$S\subseteq \bigcup^m_{i=1}U(h,t_i):=W.$$ 
 Let $\bar{\tau}=\max_{1\leq i\leq m} t_i$ and take $T=\{t_1,\ldots,t_m\}$. For any $\delta>0$, set $W_\delta=\{x:V(x)<\delta\}$, and choose $p>0$ such that 
 $S\subseteq W_p\subseteq W$. 
 Let $\mathcal{U}=W_p$.
Again using the arguments similar to the proof of \cite[Theorem 1]{jiang2015}, one  can prove that 
for all $x\in \mathcal{U}\setminus S$, there is a sequence of real numbers $\tau_i\uparrow +\infty$ with $\tau_0=0$ and $\tau_{i+1}-\tau_i\in T$, $i=0,1,2\ldots$, such that
\begin{equation}\label{equ:V_tau_i}
    V(\varPsi_{\tau_{i+1}}(x))<(1-h)V(\varPsi_{\tau_{i}}(x))<\cdots<(1-h)^{i+1}V(x),
\end{equation}
and for all $t\geq 0$, there exist $\tau_i$ and $\tau_{i+1}$ with $\tau_i\leq t< \tau_{i+1}$, such that
$V(\varPsi_t(x))\leq\bar{\beta} V(\varPsi_{\tau_{i}}(x))$, where
$$\bar{\beta}:=\displaystyle\sup_{0\leq t\leq\bar{\tau}}\sup_{x\in M}\zeta_M(t,x)+2\leq C^{l_{\bar{\tau}}+1}+2<+\infty.$$ 
Therefore, $V(\varPsi_t(x))\rightarrow 0$ as $t\to +\infty$ by \eqref{equ:V_tau_i} for all $x\in \mathcal{U}\setminus S$. It follows that $\omega(x)\subseteq S$ for all $x\in \mathcal{U}$, that is, $S$ attracts all points of $\mathcal{U}$.

On the other hand, for any neighborhood $\mathcal{V}$ of $S$, there exists a $\mu>0$ such that $\bar{\beta}\mu<p$ and $W_{\bar{\beta}\mu}\subseteq \mathcal{V}$, that is, $W_{\bar{\beta}\mu}\subseteq W_p\cap \mathcal{V}\subseteq \mathcal{U}$. Set $\mathcal{U}_1=W_\mu$. Note that $\mathcal{U}_1\subseteq W_{\bar{\beta}\mu}\subseteq \mathcal{U}$ because $\bar{\beta}>1$. Therefore, for all $x\in \mathcal{U}_1\setminus S$, one has 
$$V(\varPsi_t(x))\leq\bar{\beta} V(\varPsi_{\tau_{i}}(x))\leq\bar{\beta} V(x)<\bar{\beta}\mu, ~~t\geq 0,$$
and it follows that $\varPsi_t(x) \in W_{\bar{\beta}\mu}$ for all $t\geq 0$. Thus, $\varPsi_t(\mathcal{U}_1)\subseteq \mathcal{V}$ for all $t\geq 0$, that is, $S$ is stable under $\varPsi$. Then by \cite[Lemma 3.3.1]{hale1988}, there is a neighborhood $\mathcal{W}_S$ of $S$ in $M$ such that $S$ attracts every compact set $\mathcal{K
}$ of $\mathcal{W}_S$ (though  \cite[Lemma 3.3.1]{hale1988} deals with invariant sets, the
same proof works for positively invariant sets), and hence, $\omega(\mathcal{K})\subset S$. By the positive invariance of $S$, one has $$J=\bigcap_{t\geq 0}\varPsi_t(S)=\omega(S)$$
is an invariant compact set with $\omega(\mathcal{K})\subset J$ for all compact set $\mathcal{K
}\subset\mathcal{W}_S$. Therefore, $J\subset S$ attracts the compact neighborhood of $S$ contained in $\mathcal{W}_S$, and hence $J$ is an attractor.
\end{proof}

\section{Permanence for Kolmogorov systems}\label{sec:4}
We now turn to the case when the dynamical system is generated by the Kolmogorov system of differential equations
\begin{equation}\label{ecological-sys}
\dot{x}_i=x_if_i(x),\quad x_i\geq 0,\quad i=1,\ldots,n
\end{equation}
in the nonnegative cone $\mathbb{R}_+^n$, 
where $f_i:\mathbb{R}_+^n \to \mathbb{R}$ are $C^1$. Here, $x_i$ denotes the population density and $f_i$ the per capita growth rate of species $i$. To avoid technical conditions assume that $f=(f_1,\ldots,f_n)$ is such that solutions to system \eqref{ecological-sys} are defined for all $t\geq 0$ and continuous with respect to the time and initial values. Let $\varPsi:\mathbb{R}_+\times\mathbb{R}_+^n \to \mathbb{R}_+^n$ denote the resulting semi-flow that describes the evolution of states in positive time, where the solution of \eqref{ecological-sys} with initial value $x\geq 0$ is given by $x(t) = \varPsi_t(x)$. Note that $0$ is a trivial equilibrium for $\varPsi$.
\begin{definition}[\cite{hofbauer1987coexistence,Hutson1984,Hutson1982}]
The semi-flow $\varPsi$ generated by system \eqref{ecological-sys} {\rm(}or system \eqref{ecological-sys}{\rm)} is said to be permanent if there exists a compact positively invariant set $K\subseteq \dot{\mathbb{R}}^n_+$ such that for every $x\in \dot{\mathbb{R}}^n_+$, there exists a $t_0=t_0(x)>0$ such that the tail $\gamma^+_{t_0}(x)\subseteq K$. If system \eqref{ecological-sys} is not permanent, then it is said to be impermanent.
\end{definition}
\begin{remark}
Since $K\subseteq \dot{\mathbb{R}}^n_+$ is compact, the distance $d(K,\partial \mathbb{R}^n_+)$ of $K$ from the boundary $\partial \mathbb{R}^n_+$ is thus non-zero, and $\omega(x)\subseteq K$ for all $x\in \dot{\mathbb{R}}^n_+$. Equivalently, permanence means that there exist $\delta,D>0$ such that 
\begin{equation}
\delta\leq \liminf_{t\to +\infty}\varPsi_i(t,x)\leq \limsup_{t\to +\infty}\varPsi_i(t,x) \leq D,\quad i=1,\ldots,n,
\end{equation}
for all $x\in \dot{\mathbb{R}}^n_+$ {\rm(\cite{hofbauer1998,Kon2004,Smith2011})}. In a permanent system, species can coexist permanently in the sense that when the population densities of all species are positive, each population density will be eventually bounded away from zero and infinity. As a consequence extinction and explosion cannot occur.
\end{remark}

Consider a $C^1$ function $V: \mathbb{R}_+^n \to \mathbb{R}$. Use the notation $\dot{V}(x)$ to denote the total derivative of ${V}(x)$ along the vector field defined by \eqref{ecological-sys}, that is, 
$$
\dot{V}(x):=\sum_{i=1}^nx_if_i(x)\frac{\partial V(x)}{\partial x_i}.
$$

Set $V(x)=x^{\nu_1}_1x^{\nu_2}_2\cdots x^{\nu_n}_n$ for some $\nu_1,\ldots,\nu_n>0$. Then 
$$
\dot{V}(x)=V(x)\sum_{i=1}^n\nu_if_i(x).
$$
Note that $V(x) = 0$ for all $x\in \partial \mathbb{R}_+^n$ and $V(x) > 0$ for all $x \in\dot{\mathbb{R}}_+^n$. 
For $x\in \dot{\mathbb{R}}_+^n$, let
\begin{equation}\label{equ-gx}
g(x):=\frac{\dot{V}(x)}{V(x)}=\sum^n_{i=1}\nu_if_i(x).
\end{equation}
Since $f_i:\mathbb{R}^n_+\to \mathbb{R}$ are continuous, \eqref{equ-gx}
provides a continuous extension of $g$ to $\mathbb{R}^n_+$, that is
\begin{equation}\label{gx}
g(x)=\sum^n_{i=1}\nu_if_i(x), \quad x\in \mathbb{R}^n_+, 
\end{equation}
and $g$ is continuous on $\mathbb{R}^n_+$.

\subsection{Permanence for dissipative systems} 
We now present the criteria on permanence and impermanence for dissipative system \eqref{ecological-sys} with a global attractor $J$. For dissipative system \eqref{ecological-sys}, the solutions in $\mathbb{R}_+^n$ are defined for all $t\geq 0$, and hence it generates a dissipative semi-flow $\varPsi:\mathbb{R}_+\times\mathbb{R}_+^n \to \mathbb{R}_+^n$. 
\begin{theorem}\label{theorem:permanence}
Suppose that the semi-flow $\varPsi$ generated by \eqref{ecological-sys} is dissipative with a global attractor $J$. Set $\partial J=J\cap \partial\mathbb{R}_+^n$. If there are real numbers $\nu_1,\ldots,\nu_n>0$ such that
\begin{equation}\label{con-01}
  \sup_{t\geq 0} \int_0^t g(\varPsi_s(x)) ds>0,\quad\forall x\in\Lambda(\partial J),
\end{equation}
where $g$ is given by \eqref{gx}, then $\varPsi$ is permanent; if instead 
\begin{equation}\label{con-02}
   \inf_{t\geq 0} \int_0^t g(\varPsi_s(x)) ds<0,\quad\forall x\in\Lambda(\partial J),
\end{equation}
then $\varPsi$ is impermanent.
\end{theorem}
\begin{proof}

Let 
\begin{equation}\label{fun-Vx}
V(x)=x^{\nu_1}_1x^{\nu_2}_2\cdots x^{\nu_n}_n, \quad x\in \mathbb{R}^n_+.
\end{equation}
Note that $V(x)=0$ if and only if $x\in \partial \mathbb{R}^n_+$. We first show that if \eqref{con-01} holds then $\varPsi$ is permanent. For $x\in \dot{\mathbb{R}}^n_+$, we let $\theta(t,x)=V(\varPsi_t(x))/V(x)$, $t\geq 0$.
Then 
\begin{equation}\label{theta1}
\theta(t,x)=\exp\left\{ \int_0^t \frac{\dot{V}(\varPsi_s(x))}{V(\varPsi_s(x))} ds\right\}=\exp\left\{\int_0^t g(\varPsi_s(x))ds\right\}, ~x\in \dot{\mathbb{R}}^n_+.
\end{equation}
Since $g\circ \varPsi:\mathbb{R}_+\times \mathbb{R}^n_+\to \mathbb{R}$ is continuous, \eqref{theta1}
provides a continuous extension of $\theta(t,\cdot)$ to $\mathbb{R}^n_+$. 
Consider the semi-flow $\varPsi|_J:\mathbb{R}_+\times J\to J$ and $V:J\to \mathbb{R}_+$ given by \eqref{fun-Vx}, where $\varPsi|_J$ denotes the restriction of $\varPsi$ on $\mathbb{R}_+\times J$. Note that $V(x)=0$ if and only if $x\in \partial J$. For $x\in \partial J$, one has
\begin{equation}\label{equ-vartheta}
\vartheta_J(t,x):=\displaystyle\liminf_{\mbox{\tiny$\begin{array}{c}
y\rightarrow x\\
y\in J\setminus \partial J
\end{array}$}}\frac{V(\varPsi_t(y))}{V(y)} 
=\displaystyle\liminf_{\mbox{\tiny$\begin{array}{c}
y\rightarrow x\\
y\in J\setminus \partial J
\end{array}$}} \theta(t,y)=\theta(t,x).
\end{equation}
Let $\psi_J(x):=\sup_{t\geq 0}\vartheta_J(t,x)$ for $x\in \partial J$. 
It follows from \eqref{con-01} and \eqref{equ-vartheta} that
\begin{equation}\label{equ:g_phi}
    \psi_J(x)=\sup_{t\geq 0}\vartheta_J(t,x)=\sup_{t\geq 0}\theta(t,x)=\sup_{t\geq 0}\exp\left\{\int_0^t g(\varPsi_s(x))ds\right\}>1
\end{equation}
for all $x\in\Lambda(\partial J)$. Obviously, $\psi_J(x) >0$ for all $x\in \partial J$. In Lemma \ref{repellor_T} take $M=J$, $S=\partial J$ and $\varPsi=\varPsi|_J$. By Remark \ref{con:repel}, one has
$\psi_J(x)>1$ for all $x\in \partial J$.

Let $W=O_\varepsilon(J)$ be an $\varepsilon$-neighborhood of $J$ in $\mathbb{R}^n_+$. Since $J$ is the global attractor, for all $x\in \mathbb{R}^n_+$, $\gamma^+(x)\cap W\neq \emptyset$. Then by \cite[Lemma 2.1]{Hutson1984}, $\gamma^+(\overline{W})$ is compact and positively invariant such that for any $x\in \mathbb{R}^n_+$, there exists a $\tau=\tau(x)>0$ satisfying $\gamma^+_\tau(x)\subseteq \gamma^+(\overline{W})$. Now take $M=\gamma^+(\overline{W})$ and $S=M\cap \partial \mathbb{R}^n_+$. Consider the semi-flow $\varPsi|_M:\mathbb{R}_+\times M\to M$ and the function $V:M\to \mathbb{R}_+$ given by \eqref{fun-Vx}. Clearly, $S$ and $M\setminus S$ are positively invariant under $\varPsi|_M$.  Note that $J\subseteq M, \partial J\subseteq S$ and $V(x)=0$ if and only if $x\in S$. By \eqref{equ-gx}, 
$$\dot{V}(x)=V(x)g(x)$$
for all $x\in M\setminus S$. Since $\omega(x)\subseteq J$ for all $x\in \mathbb{R}^n_+$, one has 
$$\omega(x)\subseteq J\cap \partial \mathbb{R}^n_+=\partial J,~~ \forall x\in S,$$
by the definition of \eqref{ecological-sys}. 
Thus, $\Lambda(S)\subseteq \partial J$, so by recalling that $\psi_J(x)>1$ for all $x\in \partial J$ one has
$$
\sup_{t\geq 0}\int_0^t g(\varPsi_s(x))ds>0$$
for all $x\in \Lambda(S)$. It then follows from \cite[Corollary 2]{fonda1988} that $S$ is a repellor in $M$. Therefore, $\varPsi$ is permanent by Remark \ref{con:repel} and  recalling that for any $x\in \dot{\mathbb{R}}^n_+$, there exists a $\tau=\tau(x)>0$ such that $\varPsi_\tau(x) \in M$. 

Now suppose that \eqref{con-02} holds. Consider the semi-flow $\varPsi|_J:\mathbb{R}_+\times J\to J$ and the function $V:J\to \mathbb{R}_+$ given by \eqref{fun-Vx}. In Theorem \ref{attractor_T} take $M=J,S=\partial J$ and $\varPsi=\varPsi_J$. Fix some $t_0>0$. 
Since $g\circ \varPsi:\mathbb{R}_+\times M\to \mathbb{R}$ is continuous, there exists a $\mu>0$ such that $g(\varPsi_t(x))<\mu$ for all $(t,x)\in [0,t_0]\times M$. Thus, there exists a $C>1$ such that $\theta(t,x)\leq \exp\left\{\mu t_0\right\}\leq C$
for all $0\leq t \leq t_0$ and $x\in M$, and hence $\frac{V(\varPsi_t(x))}{V(x)}\leq C$ for all $x\in M\setminus S$ and $0\leq t\leq t_0$. It follows from \eqref{con-02} that 
$$
\varphi_J(x)=\inf_{t\geq 0}\zeta_J(t,x)=\inf_{t\geq 0}\theta(t,x)=\inf_{t\geq 0}\exp\left\{\int_0^t g(\varPsi_s(x))ds\right\}<1
$$
for all $x\in\Lambda(\partial J)$, where $\zeta_J(t,x)$ is defined in Theorem \ref{attractor_T}. It follows from Theorem \ref{attractor_T} that $\partial J$ is an attractor for $\varPsi_J$. Therefore, there exists some $x\in \dot{\mathbb{R}}^n_+$ such that $\omega(x)\subseteq \partial J\subseteq \partial \mathbb{R}^n_+$, and hence $\varPsi$ is impermanent.
\end{proof}

\begin{corollary}\label{coro:permanence}
Suppose that the semi-flow $\varPsi$ generated by \eqref{ecological-sys} is dissipative with a global attractor $J$. Set $\partial J=J\cap \partial\mathbb{R}_+^n$. If there are real numbers $\nu_1,\ldots,\nu_n>0$ such that
\begin{equation}\label{con-01-1}
    g(x)=\sum^n_{i=1}\nu_if_i(x)>0,\quad\forall x\in\Lambda(\partial J),
\end{equation}
where $\partial J=J\cap \partial\mathbb{R}_+^n$, then $\varPsi$ is permanent; if instead 
\begin{equation}\label{con-02-2}
    g(x)=\sum^n_{i=1}\nu_if_i(x)<0,\quad\forall x\in\Lambda(\partial J),
\end{equation}
then $\varPsi$ is impermanent.
\end{corollary}
\begin{proof}
Since $\Lambda(\partial J)$ is positively invariant under $\varPsi$, \eqref{con-01-1} (resp. \eqref{con-02-2}) implies that \eqref{con-01} (resp. \eqref{con-02}) holds. Therefore, the conclusion follows from Theorem \ref{theorem:permanence}.
\end{proof}

\subsection{Permanence for competitive systems}
By the carrying simplex theory provided by Hirsch in \cite{Hirsch2001}, every competitive (i.e. $\frac{\partial f_i}{\partial x_j}<0$ for all $x\in \mathbb{R}_+^n$) and dissipative system \eqref{ecological-sys}, for which the origin is a repellor has a globally attracting hypersurface $\Sigma$ of codimension one, called the carrying simplex in \cite{Zeeman1994,Z993}. See \cite{hou2020} for  a weaker version of Hirsch’s result provided recently by Hou. Specifically, the {\it carrying simplex} $\Sigma$ for system \eqref{ecological-sys} is the boundary in $\mathbb{R}_+^n$ of the basin of repulsion of the origin which satisfies the following properties:
\begin{enumerate}[(P1)]
\item $\Sigma\subseteq \mathbb{R}_+^n\setminus \{0\}$ is compact and invariant;
\item for any $x\in \mathbb{R}_+^n\setminus \{0\}$, there exists some $y\in \Sigma$ such that
$\displaystyle\lim_{t\to +\infty}|\varPsi_t(x)-\varPsi_t(y)|=0$;
\item $\Sigma$ is unordered, i.e. if $x,y\in \Sigma$ such that $x_i\geq y_i$ for all $i\in N$, then $x = y$;

\item $\Sigma$ is homeomorphic to the unit simplex $\Delta^{n-1}=\{x\in \mathbb{R}_+^n:\sum_i x_i=1\}$ via radial projection.
\end{enumerate}

(P1) and (P2) imply that the long-term dynamics of $\varPsi$ is accurately reflected by that in $\Sigma$; (P3) and (P4) mean that $\Sigma$ is topologically and geometrically simple. We denote the boundary of $\Sigma$, i.e. $\Sigma\cap \partial\mathbb{R}^n_+$ by $\partial \Sigma$, and the interior of $\Sigma$, i.e. $\Sigma \setminus \partial \Sigma$ by $\dot{\Sigma}$. 

Note that each $\pi_i$ is positively invariant under $\varPsi$ and $\partial\Sigma\cap \pi_i$ is the carrying simplex of $\varPsi|_{\pi_i}$, that is $\partial\Sigma$ is composed of the carrying simplices of $\varPsi|_{\pi_i}$, $i=1,\ldots,n$. $\Sigma$ contains all non-trivial equilibria, periodic orbits and heteroclinic cycles, etc. Every vertex of $\Sigma$ is an equilibrium, where $\Sigma$ and some positive coordinate axis meet, and denote by $q_{\{i\}}=q_ie_{\{i\}}$ the equilibrium at the vertex where $\Sigma$ and $\mathbb{H}_{\{i\}}^+$ meet. For one-dimensional case, the carrying simplex is a nontrivial equilibrium.

We now present the permanence and impermanence criteria for system \eqref{ecological-sys} which admits a carrying simplex.

\begin{theorem}\label{Sigma-permanence}
Suppose that the semi-flow $\varPsi$ generated by \eqref{ecological-sys} admits a carrying simplex $\Sigma$. Set $\partial \Sigma=\Sigma\cap \partial\mathbb{R}_+^n$. If there are real numbers $\nu_1,\ldots,\nu_n>0$ such that
\begin{equation}\label{Sigma:con-01}
  \sup_{t\geq 0} \int_0^t g(\varPsi_s(x)) ds>0,~~\forall x\in\Lambda(\partial \Sigma),
\end{equation}
where $g$ is given by \eqref{gx}, then $\partial \Sigma$ is a repellor within $\Sigma$  and $\varPsi$ is permanent; if instead 
\begin{equation}\label{Sigma:con-02}
   \inf_{t\geq 0} \int_0^t g(\varPsi_s(x)) ds<0,~~\forall x\in\Lambda(\partial \Sigma),
\end{equation}
then $\partial \Sigma$ is an attractor within $\Sigma$ and $\varPsi$ is impermanent.
\end{theorem}
\begin{proof}
Since $\Sigma$ is invariant and compact such that $\omega(x)\subseteq \Sigma$ for all $x\in \mathbb{R}^n_+\setminus \{0\}$, there exists an $\varepsilon$-neighborhood $O_\varepsilon(\Sigma)\subset \mathbb{R}^n_+\setminus \{0\}$ of $\Sigma$ such that $\overline{O_\varepsilon(\Sigma)}\subseteq \mathbb{R}^n_+\setminus \{0\}$ and $\gamma^+(x)\cap O_\varepsilon(\Sigma)\neq \emptyset$ for all $x\in \mathbb{R}^n_+\setminus \{0\}$. Then it follows from \cite[Lemma 2.1]{Hutson1984} that $\gamma^+(\overline{O_\varepsilon(\Sigma)})$ is a compact positively invariant set. Set $M=\gamma^+(\overline{O_\varepsilon(\Sigma)})$ and $S=M\cap \partial \mathbb{R}^n_+$. Obviously, $M$ is a compact neighborhood of $\Sigma$ and $S$, $M\setminus S$ are positively invariant under $\varPsi$. By the property (P2) of $\Sigma$, one has $\omega(x)\subseteq \Lambda(\partial \Sigma)$ for any $x\in S$, and hence $\Lambda(S)\subseteq\Lambda(\partial \Sigma)$. Suppose that \eqref{Sigma:con-01} holds. Then one has
\begin{equation}\label{S-repel}
\sup_{t\geq 0} \int_0^t g(\varPsi_s(x)) ds>0,~~\forall x\in\Lambda(S).
\end{equation} 
Now consider the semi-flow $\varPsi|_M:\mathbb{R}_+\times M\to M$. Let $V(x)=x^{\nu_1}_1\cdots x^{\nu_n}_n$, $x\in M$. Note that $V(x)=0$ if and only if $x\in S$.  By \eqref{equ-gx}, 
$$\dot{V}(x)=V(x)g(x)$$
for all $x\in M\setminus S$. Thus, $S$ is a repellor in $M$ by \cite[Corollary 2]{fonda1988}. Therefore, $\varPsi$ is permanent by Remark \ref{con:repel} and  recalling that for any $x\in \dot{\mathbb{R}}^n_+$, there exists a $\tau=\tau(x)>0$ such that $\varPsi_\tau(x) \in M$.
 
Now take $M=\Sigma$ and $S=\partial \Sigma$. By repeating the above arguments, we obtain that $\partial \Sigma$ is a repellor within $\Sigma$ if \eqref{Sigma:con-01} holds. While, by Theorem \ref{attractor_T} it is an attractor within $\Sigma$ if \eqref{Sigma:con-02} holds and hence $\varPsi$ is impermanent.
\end{proof}

\begin{remark}\label{rk:cs}
It deserves to point out that when system \eqref{ecological-sys} has a carrying simplex, it is more convenient to use it to study permanence. Note that the carrying simplex $\Sigma$ is a one-codimensional invariant subset so that $\partial \Sigma$ is a much smaller set than $\partial\mathbb{R}_{+}^{n}$. We only need to know information about the limit sets of orbits on $\partial \Sigma$, which is of course much smaller than the limit sets of all orbits on $\partial\mathbb{R}_{+}^{n}$. Moreover, as shown in the proof of Theorem \ref{Sigma-permanence}, the condition \eqref{Sigma:con-01} is sufficient to determine the boundary of the absorbing set  $\gamma^+(\overline{O_\varepsilon(\Sigma)})$ for $\dot{\mathbb{R}}_+^n$ being a repellor which guarantees the permanence.
\end{remark}

\begin{corollary}\label{Sigma:permanence}
Assume that $\varPsi$ generated by system \eqref{ecological-sys} admits a carrying simplex $\Sigma$. If there are real numbers $\nu_1,\ldots,\nu_n>0$ such that
\begin{equation}\label{Sigma:con-01-1}
    g(x)=\sum^n_{i=1}\nu_if_i(x)>0,\quad\forall x\in\Lambda(\partial \Sigma),
\end{equation}
then $\partial \Sigma$ is a repellor within $\Sigma$ and $\varPsi$ is permanent; if instead 
\begin{equation}\label{Sigma:con-02-2}
    g(x)=\sum^n_{i=1}\nu_if_i(x)<0,\quad\forall x\in\Lambda(\partial \Sigma),
\end{equation}
then $\partial \Sigma$ is an attractor within $\Sigma$ and $\varPsi$ is impermanent.

\end{corollary}
\begin{proof}
Note that \eqref{Sigma:con-01-1} (resp. \eqref{Sigma:con-02-2}) implies that \eqref{Sigma:con-01} (resp. \eqref{Sigma:con-02}) holds because $\Lambda(\partial \Sigma)$ is positively invariant. Then the conclusion follows from Theorem \ref{Sigma-permanence} immediately.
\end{proof}

\begin{corollary}\label{coro:2-permanence}
Let $n=2$. Suppose that system \eqref{ecological-sys} admits a carrying simplex $\Sigma$. If $f_i(q_{\{j\}})>0$ for $i\neq j$, $i,j=1,2$, where $q_{\{j\}}$ is the unique equilibrium on $\mathbb{H}^+_{\{j\}}$, then system \eqref{ecological-sys} is permanent. 
\end{corollary}
\begin{proof}
Note that $\partial \Sigma=\{q_{\{1\}},q_{\{2\}}\}$ and $f_i(q_{\{i\}})=0$ because $q_{\{i\}}$ is the equilibrium on $\mathbb{H}^+_{\{i\}}$. Therefore, there exist $\nu_1,\nu_2>0$ such that \eqref{Sigma:con-01-1} holds if $f_i(q_{\{j\}})>0$ for $i\neq j$, which implies the result by Corollary \ref{Sigma:permanence}.
\end{proof}
\begin{remark}
If $f_i(q_{\{j\}})\leq 0$ for some $i\neq j$, then the two-dimensional system \eqref{ecological-sys} might be impermanent. Consider the two-dimensional Lotka-Volterra system
  \begin{equation}\label{LV-2D}
     \left\{\begin{array}{l}
                  \frac{d x_1}{dt}=x_1(c_1-x_1-b_{12}x_2),\\
                 \noalign{\medskip}
                  \frac{d x_2}{dt}=x_2(c_2-b_{21}x_1-x_2),
 \end{array}\right.
\end{equation}
where $c_{i},b_{ij}>0$, $i,j=1,2$ with $c_2=b_{21}c_1$ and $b_{12}b_{21}<1$. It is easy to check that there is no positive equilibrium in system \eqref{LV-2D} and the equilibrium $(0,c_2)$ is unstable. Therefore, the equilibrium $(c_1,0)$ is globally attracting, and hence it is impermanent.
\end{remark}
\begin{corollary}\label{coro:3-permanence}
Let $n=3$. Suppose that system \eqref{ecological-sys} admits a carrying simplex $\Sigma$. If there are real numbers $\nu_1,\nu_2,\nu_3>0$ such that 
\begin{equation}\label{Sigma-con-02}
    g(\hat{x})=\sum^3_{i=1}\nu_if_i(\hat{x})>0\,(resp. <0),\quad\forall \hat{x}\in \mathcal{E}(\varPsi)\cap \partial \Sigma,
\end{equation}
then $\partial \Sigma$ is a repellor \,{\rm(}resp. an attractor {\rm)} within $\Sigma$ and $\varPsi$ is permanent \,{\rm(}resp. impermanent{\rm)}.
\end{corollary}
\begin{proof}
Since $\partial \Sigma\cap \pi_i$ is an invariant curve for each $i=1,2,3$, we have $\Lambda(\partial \Sigma)\subseteq \mathcal{E}(\varPsi)$. It follows that $\Lambda(\partial \Sigma)=\mathcal{E}(\varPsi)\cap\partial \Sigma
$. Therefore, the conclusion follows from Corollary \ref{Sigma:permanence} immediately.
\end{proof}

For low-dimensional systems, one remarkable phenomenon is the occurrence of heteroclinic cycles, i.e., the cyclic arrangements of saddle equilibria and heteroclinic connections (see, for example, \cite{hofbauer1994,hofbauer1998,jiang2015,may1975}).  
\begin{figure}[h]
 \begin{center}
 \includegraphics[width=0.5\textwidth]{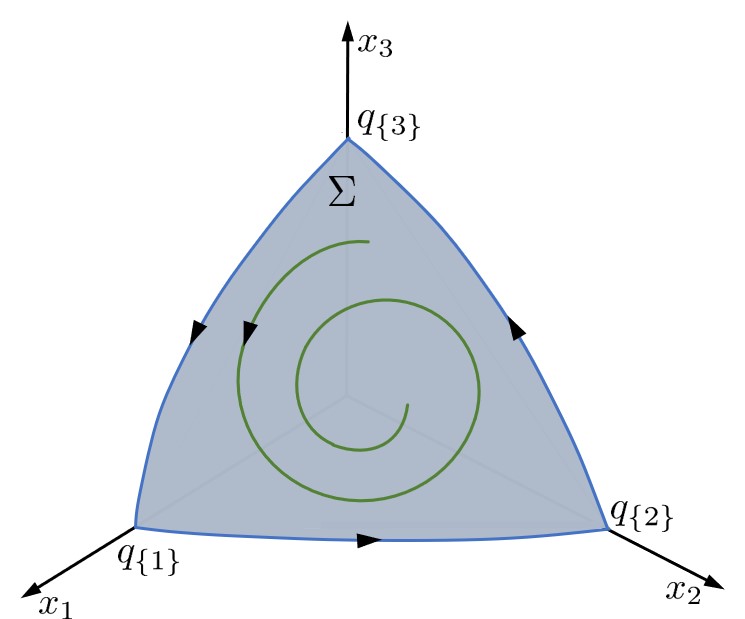}
\caption{A carrying simplex $\Sigma$ with a repelling heteroclinic cycle $\partial\Sigma$.} \label{fig:cs}
  \end{center}
\end{figure}

Let $n=3$. Suppose system \eqref{ecological-sys} admits a carrying simplex $\Sigma$ with three axial equilibria $q_{\{1\}}=(q_1,0,0)$, $q_{\{2\}}=(0,q_2,0)$ and $q_{\{3\}}=(0,0,q_3)$, which lie at the vertices of $\Sigma$. Assume that $q_{\{1\}},q_{\{2\}},q_{\{3\}}$ are saddles on $\Sigma$, and $\partial \Sigma\cap \pi_i$ is the heteroclinic connection between $q_{\{j\}}$ and $q_{\{k\}}$ (here $i,j,k$ are distinct). In this case, there are no other equilibria on $\partial \Sigma$ which is a heteroclinic cycle of May-Leonard type: $q_{\{1\}}\to q_{\{2\}} \to q_{\{3\}}\to q_{\{1\}}$ (or the arrows are reversed); see Fig. \ref{fig:cs}. This pattern means that (or all inequalities are reversed if the arrows are reversed)
\begin{subequations}\label{condition-cycle}
\begin{align}
  f_2(q_{\{1\}})\geq 0, ~f_3(q_{\{2\}})\geq 0, ~f_1(q_{\{3\}})\geq 0, \\[3pt]
  f_1(q_{\{2\}})\leq 0, ~f_3(q_{\{1\}})\leq 0, ~f_2(q_{\{3\}})\leq 0,
\end{align}
\end{subequations}
because the external eigenvalue at the axial equilibrium $q_{\{i\}}$ in direction $j$ is given by $f_j(q_{\{i\}})$, where $i\neq j$. 
Specifically, consider $q_{\{1\}}$. The Jacobian matrix at $q_{\{1\}}$ takes the form
$$
    \left[
    \begin{array}{ccc}
         x_1\frac{\partial f_1}{\partial x_1} & x_1\frac{\partial f_1}{\partial x_2} & x_1\frac{\partial f_1}{\partial x_3} \\
        \noalign{\medskip}
        0 & f_2 & 0 \\
        \noalign{\medskip}
        0 & 0 & f_3
    \end{array}
    \right],
$$
where everything is evaluated at the equilibrium. It follows that $f_2(q_{\{1\}})$ and $f_3(q_{\{1\}})$ are the external eigenvalues in directions $2$ and $3$ respectively.  It this case, $\partial \Sigma \cap \pi_3\setminus \{q_{\{2\}}\}$ is the unstable manifold of $q_{\{1\}}$ on $\Sigma$ which means $f_2(q_{\{1\}})\geq 0$, while $\partial \Sigma \cap \pi_2\setminus \{q_{\{3\}}\}$ is the stable manifold of $q_{\{1\}}$ on $\Sigma$ which means $f_3(q_{\{1\}})\leq 0$.

By applying Corollary \ref{coro:3-permanence}, we can obtain the stability of the heteroclinic cycle $\partial \Sigma$ above, which is a special case of the results provided in \cite{hofbauer1994}.

\begin{corollary}\label{coro:h_cycle}
Suppose system \eqref{ecological-sys} admits a carrying simplex $\Sigma$, and $\partial \Sigma$ is a heteroclinice cycle of May-Leonard type above.

\noindent {\rm (1)}~If
\begin{equation}\label{con:03}
    \prod_{i=1}^3 f_i(q_{\{i-1\}})+\prod_{i=1}^3f_i(q_{\{{i+1}\}})>0,
\end{equation}
 where $i\in \{1,2,3\}$ is considered cyclic, then the heteroclinic cycle $\partial \Sigma$ is a repellor within $\Sigma$, and system \eqref{ecological-sys} is permanent.

\noindent {\rm (2)}~If
\begin{equation}\label{con:04}
    \prod_{i=1}^3 f_i(q_{\{i-1\}})+\prod_{i=1}^3f_i(q_{\{{i+1}\}})<0,
\end{equation}
 where $i\in \{1,2,3\}$ is considered cyclic, then the heteroclinic cycle $\partial \Sigma$ is an attractor within $\Sigma$, and system \eqref{ecological-sys} is impermanent.
\end{corollary}
\begin{proof}
Without loss of generality, we assume that the heteroclinic cycle $\partial \Sigma$ is the type: $q_{\{1\}}\to q_{\{2\}} \to q_{\{3\}}\to q_{\{1\}}$ (see Fig. \ref{fig:cs}), and hence \eqref{condition-cycle} holds. Note that in this case $\Lambda(\partial\Sigma)=\{q_{\{1\}},q_{\{2\}},q_{\{3\}}\}$. Let $\Theta=(\Theta_{ij})_{3\times 3}$ be the characteristic
matrix of the heteroclinic cycle $\partial \Sigma$, where $\Theta_{ij} = f_j(q_{\{i\}})$, $i,j=1,2,3$. Note that $\Theta_{ii} =f_i(q_{\{i\}})=0$ because $q_{\{i\}}$ are equilibria.

For (1), we first show that there exists a positive vector $\nu=(\nu_1,\nu_2,\nu_3)^\tau$ such that
\begin{equation}\label{inequ:1}
    \sum^3_{j=1}\nu_jf_j(q_{\{i\}})>0,\quad i=1,2,3.
\end{equation} 
It is easy to see that  \eqref{inequ:1} is equivalent to $\Theta \nu\gg 0$. It follows from \eqref{condition-cycle} and \eqref{con:03} that 
$\prod_{i=1}^3 f_i(q_{\{i-1\}})>0$, and hence $f_i(q_{\{i-1\}})>0$, where $i\in \{1,2,3\}$ is considered cyclic. Let $
     \sigma=
    \left(\begin {array}{ccc}
    0&0&1\\
    1&0&0\\
    0&1&0
    \end {array} \right)
$
and $\Theta_1=\Theta\sigma$. Condition \eqref{con:03} implies that $\det\Theta>0$, and hence $\det\Theta_1>0$. Therefore, $\Theta_1$ is a matrix with positive entries occurring only in the main diagonal and its all leading principal minors are also positive. It follows that $\Theta_1$ is an $M$-matrix which implies that there exists a vector $\bar{\nu}\gg 0$, such that $\Theta_1\bar{\nu}\gg 0$. Now we can choose $\nu=\sigma \bar{\nu}$ to be the positive vector. Therefore, the heteroclinic cycle $\partial \Sigma$ is a repellor within $\Sigma$ and system \eqref{ecological-sys} is permanent by Corollary  \ref{Sigma:permanence}.

By the similar arguments as above, one can find $\nu_i>0$  such that the inequalities in \eqref{inequ:1} are reserved if \eqref{con:04} holds. Then the conclusion of (2) follows from Corollary  \ref{Sigma:permanence} immediately.
\end{proof}
\begin{remark}
We emphasize that using the carrying simplex to study the stability of the heteroclinic cycle and the permanence is equally suitable alternative to the well-known methods of \cite{hofbauer1994} and \cite[Section 16.3]{hofbauer1998}. In fact, as shown in the proof of Theorem \ref{Sigma-permanence}, the absorbing set $M=\gamma^+(\overline{O_\varepsilon(\Sigma)})$ for $\dot{\mathbb{R}}_+^3$ is a compact neighborhood of $\Sigma$. Moreover, $S=M\cap \partial \mathbb{R}^3_+$ and $M\setminus S$ are positively invariant under $\varPsi$. By the property (P2) of $\Sigma$, one has $\omega(x)\subseteq \Lambda(\partial \Sigma)$ for any $x\in S$, so $\Lambda(S)\subseteq\Lambda(\partial \Sigma)$. 
Then by \cite[Corollary 2]{fonda1988} and Theorem \ref{attractor_T}, one has that $S$ repels (resp. attracts) in $M$ if condition \eqref{con:03} (resp. condition \eqref{con:04}) holds. Therefore, if condition \eqref{con:03} holds, then there exists a neighborhood $V$ of $\partial \Sigma$ in $M$ $($hence in $\mathbb{R}_{+}^{3})$ such that given any $x \in \dot{\mathbb{R}}_{+}^{3}$ there is a $t_{0}(x)$ such that $\varPsi_t(x) \in V^{c}$ for $t \geq t_{0}(x)$, that is, $\partial \Sigma$ is repelling in $\dot{\mathbb{R}}_{+}^{3}$ and the system is permanent. While if \eqref{con:04} holds, then  there is a neighborhood $\tilde{V}$ of $\partial \Sigma$ in $\mathbb{R}_{+}^{3}$ such that $\omega(x) \subseteq \partial \Sigma$ for every $x \in \tilde{V}$, that is $\partial \Sigma$ is attracting in $\dot{\mathbb{R}}_{+}^{3}$ and the system is impermanent. Therefore, it is sufficient to study the local stability of the heteroclinic cycle within the carrying simplex to determine whether the system is permanent or not. 
\end{remark}

\section{Applications to ecological systems}
\label{sec:5}
It is the purpose of
this section to show how the criteria above enable permanent coexistence to be established for the competitive Kolmogorov system 
\begin{equation}\label{3d-equ:cssn}
    \frac{dx_i}{dt}=x_if_i(c_i,\sum_{j=1}^3b_{ij}x_j),~ c_i,b_{ij}>0,~ i,j=1,2,3
\end{equation}
 with $f_i$ satisfying \eqref{property:f} given by population models which have linearly determined nullclines. 
Note that the nullclines of system \eqref{3d-equ:cssn} given by
$$
    \frac{d x_i}{dt}=0 \Leftrightarrow x_if_i(c_i,(Bx^\tau)_i)=0 \Leftrightarrow (Bx^\tau)_i=c_i~\textrm{or}
    ~ x_i=0,~i=1,2,3.
$$
Some classical models include the  
Lotka-Volterra system (\cite{Zeeman1994,Z993})
\begin{equation}\label{L-V}
    \frac{d x_i}{dt}=x_i(c_i- \sum_{j=1}^3b_{ij}x_j),\quad i=1,2,3,
\end{equation}
Gompertz system (\cite{jnz,Lu2017})
\begin{equation}\label{Gompertz}
    \frac{d x_i}{dt}=x_i\ln\frac{c_i}{\sum_{j=1}^3b_{ij}x_j},\quad i=1,2,3,
\end{equation}
Leslie-Gower system (\cite{hou2019,JN2017})
 \begin{equation}\label{equs:LG}
 \frac{dx_i}{dt}=x_i\frac{1+c_i}{1+\sum_{j=1}^3b_{ij}x_j}-x_i,\quad i=1,2,3,
 \end{equation}
and Ricker system (\cite{hou2019,JN2017})
\begin{equation}\label{equs:Ricker}
 \frac{dx_i}{dt}=x_i\exp(c_i-\sum_{j=1}^3b_{ij}x_j)-x_i, \quad i=1,2,3.
\end{equation}

By applying Hirsch's carrying simplex theory \cite{hirsch1988}, Jiang and Niu \cite{JN2017} show that every system \eqref{3d-equ:cssn} admits a carrying simplex.
Via the existence of a carrying simplex, Jiang and Niu \cite{JN2017} generalized Zeeman's classification  in \cite{Z993} by nullcline equivalence for three-dimensional competitive Lotka-Volterra system \eqref{L-V} to any competitive systems with linearly determined nullclines given by \eqref{3d-equ:cssn}.

Let $\gamma_{ij}:=b_{ii}c_{j}-b_{ji}c_i~ (i\neq j)$ and $\beta_{ij}:=b_{ii}b_{jj}-b_{ij}b_{ji}~ (i<j)$. 
The nullcline configuration of system \eqref{3d-equ:cssn} is given by the values of $\mathrm{sgn}(\gamma_{ij})$ for distinct $i,j$, and $\mathrm{sgn}(c_k\beta_{ij}-b_{ki}\gamma_{ji}-b_{kj}\gamma_{ij})$ (if $v_{\{k\}}$ exists), modulo permutation of the indices $\{1,2,3\}$. Two systems given by \eqref{3d-equ:cssn} are said to be \textit{nullcline equivalent} if they have the same nullcline configurations. The system is \textit{nullcline stable} if and only if $\gamma_{ij}\neq 0$ and $c_k\beta_{ij}-b_{ki}\gamma_{ji}-b_{kj}\gamma_{ij}\neq 0$ {\rm(}if $v_{\{k\}}$ exists{\rm)}. Note that if there is a unique planar equilibrium 
$v_{\{k\}}$ in the interior of coordinate plane $\pi_k$, then $\beta_{ij}\neq 0$, and moreover, $\beta_{ij}> 0$  (resp. $\beta_{ij}<0$) if $v_{\{k\}}$ attracts (resp. repels) on $\partial \Sigma\cap \pi_k$, where $i,j,k$ are distinct. An equivalence class is said to be \textit{stable} if each system in it is nullcline stable. Specifically, Jiang and Niu proved in \cite{JN2017} that there is a total of 33 stable nullcline classes for system \eqref{3d-equ:cssn}, which are described in terms of inequalities on parameters $b_{ij},c_i$, and independent of $f_i$; see \cite[Appendix A]{JN2017}. In particular, the boundary of the carrying simplex in class 27 is a heteroclinic cycle of May-Leonard type.

Our subject here is to show how to obtain the parameter conditions on $b_{ij}$ and $c_i$ guaranteeing the permanence or impermanence for systems of type \eqref{3d-equ:cssn} by the criteria established in Section \ref{sec:4}. 

\begin{figure}[h]
 \begin{center}
 \includegraphics[width=0.32\textwidth]{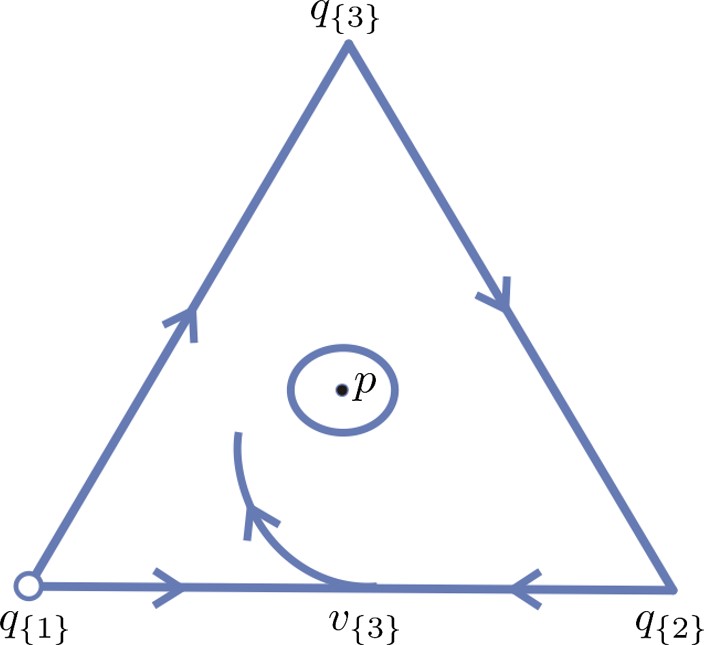}
\caption{The dynamics on $\Sigma$ for class $29$. An equilibrium is represented by an open dot $\circ$ if it is a repellor on $\Sigma$, while if the equilibrium is a saddle on $\Sigma$, we do not use any symbol.} \label{fig:c29}
  \end{center}
\end{figure}

\begin{lemma}[\cite{JN2017}]\label{class-29}
A system \eqref{3d-equ:cssn} is in the class $29$ if and only if
\begin{itemize}
\item[{\rm(i)}] $\gamma_{12}>0, \gamma_{13}>0, \gamma_{21}>0, \gamma_{23}<0, \gamma_{31}<0, \gamma_{32}>0$;
\item[{\rm(ii)}] $b_{31}\gamma_{21}+b_{32}\gamma_{12}-c_3\beta_{12}<0$.
\end{itemize}
The phase portrait on the carrying simplex $\Sigma$ for class $29$ is as shown in Fig. \ref{fig:c29}. 
\end{lemma}

\begin{proposition}\label{permanence-29}
Any system \eqref{3d-equ:cssn} in the class $29$ is permanent.
\end{proposition}
\begin{proof}
Let $f_i(x)$ denote $f_i(c_i,\sum_{j=1}^3b_{ij}x_j)$ in \eqref{3d-equ:cssn}, $i=1,2,3$. For the system in class $29$, there are four equilibria on $\partial \Sigma$, that is, $q_{\{1\}}$, $q_{\{2\}}$, $q_{\{3\}}$ and $v_{\{3\}}$; see Fig. \ref{fig:c29}. 
By Corollary \ref{coro:3-permanence}, it suffices to prove that there are real numbers $\nu_1,\nu_2,\nu_3>0$ such that the following inequalities hold:
\begin{subequations}\label{inequalities}
\begin{align}
  \nu_1f_1(q_{\{1\}})+\nu_2f_2(q_{\{1\}})+\nu_3f_3(q_{\{1\}})>0;\label{ineq:aa} \\[2pt]
  \nu_1f_1(q_{\{2\}})+\nu_2f_2(q_{\{2\}})+\nu_3f_3(q_{\{2\}})>0;\label{ineq:bb}\\[2pt]
  \nu_1f_1(q_{\{3\}})+\nu_2f_2(q_{\{3\}})+\nu_3f_3(q_{\{3\}})>0;\label{ineq:cc}\\[2pt]
  \nu_1f_1(v_{\{3\}})+\nu_2f_2(v_{\{3\}})+\nu_3f_3(v_{\{3\}})>0.
  \label{ineq:dd}
\end{align}
\end{subequations}
For an equilibrium $\hat{x}$ of system \eqref{3d-equ:cssn}, one has $f_i(\hat{x})=0$ for all $i\in \kappa(\hat{x})$. Therefore, 
$f_i(q_{\{i\}})=0$, $i=1,2,3$ and $f_1(v_{\{3\}})=f_2(v_{\{3\}})=0$. Note that $v_{\{3\}}$ attracts on $\partial \Sigma\cap \pi_3$ (see Fig. \ref{fig:c29}), so $\beta_{12}>0$. Then by condition (ii) in Lemma \ref{class-29}, one has
$$ b_{31}\gamma_{21}+b_{32}\gamma_{12}-c_3\beta_{12}<0 \Leftrightarrow  f_3(v_{\{3\}})>0.$$ 
So, \eqref{ineq:dd} holds for any $\nu_1,\nu_2,\nu_3>0$.
Since $\gamma_{12},\gamma_{13}>0$ by condition (i) in Lemma \ref{class-29}, one has $f_2(q_{\{1\}}),f_3(q_{\{1\}})>0$. Thus, \eqref{ineq:aa} holds for any $\nu_1,\nu_2,\nu_3>0$. The inequalities \eqref{ineq:bb} and \eqref{ineq:cc} can be written as 
\begin{subequations}\label{inequalities-2}
\begin{align}
  \nu_1f_1(q_{\{2\}})+\nu_3f_3(q_{\{2\}})>0;\label{ineq:aa-2}\\[2pt]
  \nu_1f_1(q_{\{3\}})+\nu_2f_2(q_{\{3\}})>0.\label{ineq:bb-2}
\end{align}
\end{subequations}
We first fix a $\nu_2>0$. It follows from $\gamma_{32}>0$ that $f_2(q_{\{3\}})>0$, so for sufficiently small $\nu_1>0$ one has \eqref{ineq:bb-2} holds. Now fix some $\nu_1>0$ such that \eqref{ineq:bb-2} holds. Note that $\gamma_{21}>0$, so $f_1(q_{\{2\}})>0$. Then we can choose a $\nu_3>0$ sufficiently small such that \eqref{ineq:aa-2} holds.  Such $\nu_1,\nu_2,\nu_3>0$ ensure that the inequalities \eqref{ineq:aa}--\eqref{ineq:dd} hold. This proves that each system in class $29$ is permanent. 
\end{proof}

\begin{proposition}\label{permanence-31-33}
Any system \eqref{3d-equ:cssn} in the classes $31$ and $33$ is permanent.
\end{proposition}
\begin{proof}
The proof is similar to that of Proposition \ref{permanence-29}, so we omit it here.
\end{proof}

\begin{proposition}\label{permanence-27}
Assume that system \eqref{3d-equ:cssn} is in class $27$. Let 
\begin{equation}
\varrho=\lambda_{12}\lambda_{23}\lambda_{31}+\lambda_{21}\lambda_{13}\lambda_{32}
\end{equation}
with $\lambda_{ij}=f_j(c_j,(Bq^\tau_{\{i\}})_j)$, $i\neq j$. 
Then the system is permanent if $\varrho>0$, while it is impermanent if $\varrho<0$.
\end{proposition}
\begin{proof}
Since the boundary of carrying simplex is a heteroclinic cycle for system \eqref{3d-equ:cssn} in class $27$, the conclusion follows from Corollary \ref{coro:h_cycle}.
\end{proof}

\begin{theorem}\label{permanence-3D}
System \eqref{3d-equ:cssn} is permanent if it is in classes 29, 31, 33 and class 27 with $\varrho>0$, while it is impermanent if it is in classes 1--26, 28, 30, 32 and class 27 with $\varrho<0$.
\end{theorem}
\begin{proof}
For each system in classes 1--26, $28$, $30$ and $32$, there always exists an equilibrium on $\partial \Sigma$ which is an attractor on $\Sigma$, so it is impermanent. The rest results follow from Propositions \ref{permanence-29}--\ref{permanence-27}. See \cite[Appendix A]{JN2017} for the parameter conditions of each class. 
\end{proof}

\section{Discussion}
In this paper, we study the permanence and impermanence for the competitive Kolmogorov systems \eqref{kolmogorov-equas}. We first provide a criterion on the existence of attractors via the average Liapunov function in continuous-time dynamical systems, which is an extension of Hutson's result \cite{Hutson1984} on the existence of repellors in continuous-time dynamical systems. The result is particularly easy to apply to establish impermanence criteria for Kolmogorov systems \eqref{kolmogorov-equas}, which is a continuous version of the one provided in \cite{jiang2015} for discrete-time dynamical systems. 

We then give a general result for the stability of the boundary of the carrying simplex for continuous-time competitive Kolmogorov systems \eqref{kolmogorov-equas}, which determines the permanence and impermanence of the systems. Specifically, the boundary of the carrying simplex being a repellor implies that the system is permanent, while the existence of an attractor on the boundary of the carrying simplex implies that the system is impermanent. The conditions of our criteria are very easy to verify, and in particular, they are finitely computable conditions which only depend on the nontrivial boundary equilibria for three-dimensional competitive Kolmogorov systems. When the boundary of the carrying simplex is a heteroclinic cycle, the criteria are just the stability criteria of the heteroclinic cycle. Therefore, in this case the system is permanent if the heteroclinic cycle is a repellor, while the system is impermanent if it is an attractor. 

Using our results, we provide a complete classification of the permanence and impermanence in terms of inequalities on parameters for the three-dimensional competitive Kolmogorov system \eqref{3d-equ:cssn} with linearly determined nullclines given by population models, which includes the classical Lotka-Volterrra system \eqref{L-V}, Gompertz system \eqref{Gompertz}, Leslie-Gower system \eqref{equs:LG} and Ricker system \eqref{equs:Ricker}. In fact, applying our criteria to each of the 33 equivalence classes for \eqref{3d-equ:cssn} given in Jiang and Niu \cite{JN2017}, we obtain that systems in their classes 29, 31, 33 and class 27 with a repelling heteroclinic cycle are permanent, while systems in their  classes 1--26, 28, 30, 32 and class 27 with an attracting heteroclinic cycle are impermanent. In particular, our classification implies that the system \eqref{3d-equ:cssn} is permanent when all the boundary equilibria are unstable if the boundary of the carrying simplex is not a heteroclinic cycle. 

We emphasize that for the system \eqref{3d-equ:cssn}, the permanence generally does not imply the global asymptotic stability of the positive equilibrium. Note that there might be limit cycles in the permanent classes for some classical models. For example, limit cycles can occur in classes 29, 31 and class 27 with a repelling heteroclinic cycle for Lotka-Volterra system \eqref{L-V} (see, for example, \cite{Gyllenberg09number, Jiang2021, Li2022, lu2002two, lu2003three, Yu2016, Z993}). While, for Gompertz system \eqref{Gompertz}, Lu and Lu \cite{Lu2017} proved that there is no limit cycle in classes 29 and 31, and consequently, the unique positive equilibrium in the two classes is globally asymptotically stable. It is an interesting problem to provide conditions under which the permanence can imply the global asymptotic stability of the positive equilibrium in these classical models, which is left for future research.

\section*{Acknowledgement}
The work was supported by the National Natural Science Foundation of China (No. 12001096) and the Fundamental Research Funds for the Central Universities (No. 2232023A-02 and 2232023G-13).

\bibliographystyle{elsarticle-num-names}
\bibliography{refs}

\end{document}